\documentclass[11pt,preprintnumbers,showpacs,preprint,showpacs,notitlepage,a4paper]{article}
\usepackage{amsmath}
\usepackage{amssymb}
\usepackage{braket}
\usepackage{graphicx}
\usepackage{dcolumn}
\usepackage{bm}
\usepackage[usenames]{color}
\usepackage{bbm,amssymb,amsmath}

\newenvironment{Claim}[1]{\vspace*{2mm}\noindent{\textit Claim} $\bf{#1}$:\hspace*{1em}}{$\Box $}
\newenvironment{proofof}[1]{\vspace*{2mm}\noindent{\textit Proof} of $\bf{#1}$:\hspace*{1em}}{$\Box $}
\newcommand{\mcR}{\mathbbm{R}}
\newcommand{\mcN}{\mathbbm{N}}
\newcommand{\mcI}{\mathbbm{I}}

\newtheorem{theorem}{Theorem}
\newtheorem{lemma}{Lemma}

\begin{document}
\title{Evaluation of Riemann Zeta function on the Line $\Re(s)=1$ and Odd Arguments}
\author{\bf{Srinivasan.A}$^{\dag}$ }

\maketitle

\begin{abstract}

\noindent  We have looked at the evaluation of the riemann zeta function at odd arguments and have provided a simple formula to approximate the value with exponential convergence. We have compared it with various other formulae present in literature. We have also evaluated an expression for the zeta function on the plane $\Re(s)=1$. 
\end{abstract}
\let\thefootnote\relax\footnotetext {$\dag$ Srinivasan.A (srinivasan1390[at]gmail[dot]com) is a final year undergraduate Engineering student at the National Institute of Technology (NIT), Karnataka(India).}

\section{Introduction} \label{intro}
The Riemann Zeta function is a widely talked about problem which features as one of the Seven Millennium Problems. Simply Stated Riemann conjectured that \textbf{'All non-trivial zeros of the Zeta function have real part equal to one-half'}. If true, the hypothesis would have profound consequences on the distribution of primes in the integers \cite{P88}. The  Analytic Continuation of Dirichlet Series or the defining Riemann Zeta function is the simply given by 
\[\zeta(s) =
\begin{cases}
\sum_{n=1}^{\infty} \frac{1}{n^s} & \Re(s) \textgreater 1 \\
\frac{1}{1-2^{1-s}}\sum_{n=1}^{\infty} \frac{(-1)^{n-1}}{n^s} & \Re(s) \textgreater 1 , s\neq 1
\end{cases} \]
Its evident that $\zeta(s)$ is an infinite series which converges when the real part of $s$ is greater than 1. It can indeed be extended meromorphically to the whole complex s-plane except for a simple pole at
$s=1$ with residue 1(\cite{T86})The Riemann zeta function satisfies the functional equation
\begin{eqnarray}
    \zeta(s) = 2^s\pi^{s-1}\ \sin\left(\frac{\pi s}{2}\right)\ \Gamma(1-s)\ \zeta(1-s) \!, 
\end{eqnarray}

The general formulae well known related to the Riemann Zeta function are
\begin{eqnarray}
  \zeta(2n) &=& (-1)^{n+1}\frac{B_{2n}(2\pi)^{2n}}{2(2n)!} \\
    \zeta(-n)&=&-\frac{B_{n+1}}{n+1}
\end{eqnarray} 
where $B_{2n}$ is a Bernoulli number; hence the values of zeta function vanishes at even negative integers because $B_m$ = 0 for all odd $m$ other than 1. \textbf{No such simple expression is known for $odd$ $positive$ $integers$}.\\

	\indent The practical importance of this zeta function is in the distribution of primes. When we look at the distribution of primes between $10^i$ and $10^{1+i}$ which $i\in \mcN$. If numbers of prime numbers lying in these different intervals is counted this distribution is related pretty closely to the zeros of the Riemann Zeta function. The distribution of primes is of prime concern in areas of security, encryption, factoring, etc. For e.g. the important RSA encryption scheme relies on the difficulty of detecting prime numbers.
\section{Uniform Convergence}
Uniform convergence is, for any $\epsilon\textgreater 0$, there exists $N$ such that for the same $N$, no matter which $x$ you pick $n\textgreater N$ implies $|f_n(x)-f(x)|\textless \epsilon$\\
	\indent The fundamental difference between point wise and uniform convergence can be best visualized by seeing that uniform convergence essentially implies that for large $n$, $f_n(x)$ is \textbf{very} close to $f(x)$ as a graph everywhere, whereas point wise convergence simply says that tracking the value of a single point over time will give the right value. \\
	\indent One example of a sequence that doesn't converge uniformly is $x_n$ on [0,1]. You can see this as, if $x \textless 1$, $x_n$ converges point wise to 0, i.e. for a given point, it converges to zero. If $x=1$, $x_n$ converges to 1 obviously. But for a given $\epsilon$, say 1/4, for any $N$ we can find a point $(\frac{1}{2})^{1/N}$ so if x is that, $x^N$=1/2.\\
\indent The 3 known methods of proving uniform convergence are Abel's Uniform Convergence Test, Weierstrass M-Test, Uniform Norm Method. We shall be employing the Uniform Norm Method to prove the uniform convergence of our functions.

\subsection{Uniform Norm Method} \label{unm}
If $A\hspace{1pt} \subseteq \hspace{1pt} R$, $f$:$A$ $\rightarrow$ $R$ , $f$ is bounded. Uniform Norm on $A$ is\\ $|| f ||_A$ : sup\{$|f(x)| : x\in A $\}
\begin{theorem}
A sequence ($f_n$)  of bounded functions on $A\hspace{1pt} \subseteq \hspace{1pt} R$ converges uniformly on $A$ to $f$ if and only if $||f_n-f||_A \rightarrow 0$
\end{theorem}\newpage
\section{Zeta Function at Odd Arguments}
\subsection{Intuition}
In recent years there has been a renewed interest in representing $\zeta(2n +1)$ (n$\in \mcN$), $\mcN$ being the
set of positive integers, by means of series which converge more rapidly than the defining series in
Section \ref{intro}. These developments seem to have stemmed from the use of the familiar series representation for the known constants like $\zeta(3)$, $\zeta(5)$, $\zeta(7)$ in \cite{SGA00},\cite{V10} \\
\begin{eqnarray}
\zeta(3)&=&-\frac{4\pi^2}{7}\sum_{k=0}^{\infty} \frac{\zeta(2k)}{(2k+1)(2k+2)2^{2k}}\\
 \zeta(5)&=& 12 \sum_{n=1}^\infty \frac{1}{n^5 \sinh (\pi n)} -\frac{39}{20} \sum_{n=1}^\infty \frac{1}{n^5 (e^{2\pi n} -1)} -\frac{1}{20} \sum_{n=1}^\infty \frac{1}{n^5 (e^{2\pi n} +1)}\! \notag \\
\zeta(7)&=&\frac{19}{56700}\pi^7 -2 \sum_{n=1}^\infty \frac{1}{n^7 (e^{2\pi n} -1)}\! 
\end{eqnarray}
The main aim in the result presented below is to reduce the complexity of the huge computational formulae presented in literature. The intuition was that something simpler would exist which could at smaller values of argument might not be exact but with increasing values give better approximations to $\zeta(2s+1)$. An explicit mention has been made in the Section \ref{notes} with the other known results for computation of $\zeta(2k+1)$. In the \cite{V10} Vestas has explicitly mentioned the computational difficulties of higher values of Bernoulli numbers. Thus we have also used an earlier result \cite{A10} to avoid the need of Bernoulli numbers due to the computational factors as mentioned in \cite{FP07}. 
\subsection{Main Result}
In general it is known that 
\begin{eqnarray}
\zeta(s)=\sum_{n=1}^{\infty} \frac{1}{n^s}&=& \frac{1}{1^s}+\frac{1}{2^s}+\frac{1}{3^s}+\frac{1}{4^s}+\dots\\
&=& (1)+(\frac{1}{2^s}+\frac{1}{4^s}+\frac{1}{6^s})+(\frac{1}{3^s}+\frac{1}{9^s}+\frac{1}{27^s})+(\frac{1}{5^s}+\frac{1}{25^s}+\frac{1}{125^s})+\dots \notag\\
 &=& (1) + \frac{\zeta(s)}{2^s} + \frac{1}{3^s-1} + \frac{1}{5^s-1} + \dots \\
&=& 1+ \frac{\zeta(s)}{2^s} + \frac{1}{k^s-1}\\
\zeta(s) &=& \frac{ \sum_{p \geq 3,p\hspace{3pt}is\hspace{3pt}prime}\frac{1}{p^s-1} +1}{1-2^{-s}}
\end{eqnarray}

$Note$: We would like to highlight that in the above formulation terms are fundamentally the poweres of primes and notions that terms such as 15 or 21 which is a product of 2 primes would not repeat. Euler Product formula connecting the Primes and the Zeta function is
\begin{eqnarray}
\zeta(s)= \Pi_{p \geq 2, p\hspace{3pt}is\hspace{3pt}prime} \frac{p^s}{p^s-1}
\end{eqnarray}
\begin{Claim}1 \label{claim1}
Assuming $p$ is a set of primes, t(s)= $\sum_{p \geq 2} \frac{1}{p^s-1}$, $ \zeta(s)_{p\geq 2}= \Pi \frac{p^s}{p^s-1}$ and given the relation\\
\begin{eqnarray}\label{eqn5}
\frac{t(2s)}{\zeta(2s)} = f \frac{t({2s+1})}{\zeta(2s+1)}
\end{eqnarray}
We can give a bound on the value of f(s) and present a convergence of the form
\begin{eqnarray}
 \lim_{s \rightarrow 0} | f(s) - 2 | \leq \epsilon \hspace{5 pt} \forall \hspace{5 pt} \epsilon \textgreater 0
\end{eqnarray}
\end{Claim} \\
The above is the main crux and its pure \textbf{observational} that $'f(s)\rightarrow 2'$. From Matlab if can be realized that the constant $f$ which links both the relations averages with higher values of $s$ is  $f=2$ . Hence, initially when $s$ is small it is different from $2$ but as $s$ increases to above 20 the value of $f$ averages to $2$. Hence during computation of $\zeta(s)$ for higher values of $s$ we can use the value $f=2$ and use the formula.\\
\begin{eqnarray}
\zeta(2s+1)=2\frac{t(2s+1)}{t(2s)} \zeta(2s)
\end{eqnarray}

Using Claim 1  with the notations of t(s) and and $\zeta(s)$ remaining consistent. 
\begin{eqnarray}
\sum_{p \geq 3} \frac{1}{p^s-1} &=& \zeta (s) (1-2^{-s})-1 \\
\sum_{p \geq 2} \frac{1}{p^s-1}(=t(s))-\frac{1}{2^s-1} &=& \zeta (s) (1-2^{-s})-1 \\
t(s)&=&\zeta(s)(1-2^{-s})-1+\frac{1}{2^s-1}
\end{eqnarray}
Using the above we can rewrite Equation \ref{eqn5}
\begin{eqnarray}
\frac{t(2s)}{\zeta(2s)} =& \frac{\zeta(2s)(1-2^{-2s})-1+\frac{1}{2^{2s}-1}}{\zeta(2s)} \\
=& 1-2^{-2s} - \frac{1}{\zeta(2s)} + \frac{1}{\zeta(2s)(2^{2s}-1)} \\
\frac{t(2s+1)}{\zeta(2s+1)} =& \frac{\zeta(2s+1)(1-2^{-(2s+1)})-1+\frac{1}{2^{2s+1}-1}}{\zeta(2s+1)}\\
=& 1-2^{-(2s+1)} - \frac{1}{\zeta(2s+1)} + \frac{1}{\zeta(2s+1)(2^{2s+1}-1)} 
\end{eqnarray}
Using Claim 1 we can rewrite the above as
\begin{eqnarray*}
\frac{t(2s)}{\zeta(2s)}&=& f \frac{t(2s+1)}{\zeta(2s+1)} \\
\Big( 1-2^{-2s} - \frac{1}{\zeta(2s)} + \frac{1}{\zeta(2s)(2^{2s}-1)} \Big) &=& f \Big(1-2^{-(2s+1)} - \frac{1}{\zeta(2s+1)} + \frac{1}{\zeta(2s+1)(2^{2s+1}-1)} \Big) 
\end{eqnarray*}
\begin{eqnarray}
\zeta(2s+1)=\frac{2-2.2^s}{(2.2^s-1)\Big(\frac{\frac{2^s-1}{2^s}+\frac{1}{\zeta(2s)}\frac{2-2^s}{2^s-1}-1}{f}-1+\frac{1}{2.2^s}\Big)}
\end{eqnarray}
Using the formula for $\zeta(2s)=\frac{(-1)^{s+1} B_{2s} (2\pi)^{2s}}{2(2s)!}$ we get the following final expression for $\zeta(2s+1) \forall s\geq 1$
\begin{eqnarray}\label{eqnf}
\zeta(2s+1)&=&\frac{f. 2^{2s+1}(1-2^{2s})}{(2^{2s+1}-1)\Big(\frac{4(2s)!(1-2^{2s-1})}{(-1)^{s+1}B_{2s} \pi^{2s} (2^{2s}-1)}+\frac{2^{2s+1}(1-f)+f-2}{2}\Big)}
\end{eqnarray}

In order to avoid problems with huge Bernoulli results we could use the Result in \cite{A10} and rewrite the Equation \ref{eqnf} as
\begin{eqnarray*}
\zeta(2k)&=&-\Big(\sum_{j=0}^{k-2}(\frac{-1}{\pi^2})^{j+1}\frac{1}{(2k-2j-1)!}\zeta(2j+2)+\frac{k}{(2k+1)!}\Big)\pi^{2k}(-1)^k\\
\zeta(2k+1) &=& -f\frac{\frac{(2^{2k+1}-1) 2^{2k}}{2(1-2^{2k}) 2(2k)!}\pi^{2k}(-1)^k}{(2-2^k)-\frac{2^{2k+1}(1-f)+f-2}{2}\Big(\sum_{j=0}^{k-2}(\frac{-1}{\pi^2})^{j+1}\frac{1}{(2k-2j-1)!}\zeta(2j+2)+\frac{k}{(2k+1)!}\Big)}
\end{eqnarray*}\\

Numerical Results to authenticate the above formula (Using $f =2$) \footnote{We have used the value of $f=2$ for computing the table below. This provides us accurate results for higher values of computation of $\zeta(s)$}

\begin{tabular}{| l | c | c | c | }
\hline
 & Formula Result & Actual Result & Difference \\
$\zeta$(3)  & 1.21992 &  1.202056 & 1.7861 $\times$ $10^{-2}$ \\
$\zeta$(5)  & 1.03933 &  1.036927 & 2.3021 $\times$ $10^{-3}$ \\
$\zeta$(7)  & 1.00861 &  1.008349 & 2.4187 $\times$ $10^{-4}$ \\
$\zeta$(9)  & 1.00204 &  1.002008 & 2.5985 $\times$ $10^{-5}$ \\
$\zeta$(11) & 1.00494 & 1.004941  & 2.8476 $\times$  $10^{-6}$ \\
$\zeta$(13) & 1.00012 & 1.000122  & 3.1468 $\times$ $10^{-7}$ \\
$\zeta$(15) & 1.00003 & 1.000030  & 3.4890 $\times$ $10^{-8}$ \\
\hline
\end{tabular}\\

\subsection{Notes} \label{notes}

There are plenty formulae for Riemann zeta numbers at odd numbers but what differentiates the above from the rest of them are: \\
\begin{itemize}
\item The convergence is $Exponential$ in Nature. We can see that the difference of the actual zeta values at odd numbers and the zeta values according to the formula proposed by us. We state that with higher values of zeta the difference would be almost negligible with orders of $O(10^{-k})$. 
\item The $Simplicity$ in evaluation. A few formulae presented in literature from \cite{ZW93},\cite{SGA00},\cite{S98} are as follows :

\begin{multline}
\zeta(2m+1)=\frac{(-1)^m \pi^{2m}}{(1-2^{-2m})}\Big(-\frac{\log(2)}{\Gamma(2m+2)}  
+ \sum_{n=1}^{\infty}\frac{(2-2^{1-2n}) \Gamma(2n) \zeta(2n)}{\Gamma(2m+2n+2)} \\
+ \frac{1}{1-2^{-2m}}\sum_{n=1}^{m-1}\frac{(2^{2n-2m}-1)-(\pi^2)^n\zeta(2m-2n+1)}{\Gamma(2n+2)}\Big) n\in \mcN
\end{multline}
\begin{multline}
\zeta(2n+1)=(-1)^{n-1}\frac{(2\pi)^{2n}}{(2n)!(2^{2n+1}-1)}\Big(\log(2)+\sum_{k=0}^{\infty}\frac{\zeta(2k)}{(k+n)2^{2k}}\\
+(2n)!\sum_{j=1}^{n-1}\frac{(-1)^j}{(2n-2j)!} \Big(\frac{2^{2j}-1}{(2\pi)^{2j}}\Big)\zeta(2j+1)\Big) n\in \mcN
 \end{multline}

\begin{multline}
\zeta(2n+1)=(-1)^{n-1}\frac{(2\pi)^{2n}}{(2n)!(3^{2n+1}-1)}\Big(\log(3)+2\sum_{k=0}^{\infty}\frac{\zeta(2k)}{(k+n)3^{2k}}\\
+(2n)!\sum_{j=1}^{n-1} \frac{(-1)^j}{(2n-2j)!}\Big(\frac{3^{2j}-1}{(2\pi)^{2j}}\Big) \zeta(2j+1) \\
-\frac{(2n)!}{\sqrt(3)}\sum_{j=1}^{n}\frac{(-1)^j}{(2n-2j+1)!}\frac{2\zeta(2j,\frac{1}{3})-(3^{2j}-1)\zeta(2j)}{(2\pi)^{2j-1}}\Big)  n\in \mcN
\end{multline}

\begin{multline}
\zeta(2n+1)=(-1)^{n-1}\frac{(2\pi)^{2n}}{(2n)!(2^{4n+1}+2^{2n}-1)} \Big(\log(2)\\ +\sum_{k=0}^{\infty}\frac{\zeta(2k)}{(k+n)4^{2k}} +(2n)!\sum_{j=1}^{n-1}\frac{(-1)^j}{(2n-2j)!} 
\Big( \frac{2^{2j}-1}{(2\pi)^{2j}} \Big) \zeta(2j+1) \\
-(2n)!\sum_{j=1}^{n}\frac{(-1)^j}{(2n-2j+1)!}.  \frac{\zeta(2j,\frac{1}{4})-2^{2j-1}(2^{2j}-1)\zeta(2j)}{(2\pi)^{2j-1}}\Big) n\in \mcN
\end{multline}

Quite explicitly it can be seen that our formula is  $independent$ of the $huge$ $number$ of $double$ $sums$ and $integrals$ which occur in the other evaluations of $\zeta(2n+1)$. The convergence of the formula mentioned might be exponential and little slower than other formulae present in literature but we feel reducing complexity(by reducing the parameters, increasing ease of calculation by softwares like Mathematica, Matlab, Maple etc.) helps in the same evaluation time of higher values of zeta functions at odd arguments. We also negate the need for Bernoulli Numbers which at higher values might be huge and cause extra evaluation time. Time complexity of Bernoulli Numbers has been explicitly mentioned in \cite{FP07}. In a certain ways our expression also resembles a few of them presented above.

\item $Preknowledge$  of $Odd \hspace{4pt} \zeta$  Values isn't  Required: It can be seen in the most of the previous formulae the sums include terms which require evaluation and hence pre-knowledge of $\zeta$ at odd values. Thus our formula can be used to evaluate any $\zeta(k)$ without the need for knowing the values of $\zeta(2k-1)$ for all values of k. 
\item Exact expressions for $\zeta$ values at odd arguments are generally known for the famous $\zeta(3), \zeta(5), \zeta(7)$ where the convergence is tremendously fast and accurate. As such most of the formulae for general expressions of $\zeta(2k+1)$ are mainly of approximations. And we would like to mention that ours is one among the many approximations with exponential convergence and easier computation.
\item The presented formulation doesnt converge rapidly for higher values of zeta function purely because of the convergence of the of the zeta function to 1 with higher parameters, cause in that situation the convergence would
uniform convergence and not exponential in nature. The exponential nature of the convergence suggests that the presented formula is an alternative expression for the zeta function at odd values.
\end{itemize}

\section{Riemann Zeta function at $Re(s)=1$}

\subsection{Intuition}
In general it is known that $\zeta(s)$ diverges for $s=1$ and there are many definitions which meromorphically extend it to the entire plane except at the pole. We wanted to investigate the nature of the nature $\zeta(s)$ for values of the form $s=1+ib$. It obviously isnt an unbounded value and turns out to be bounded. There are quite a few stringent uniform continuity conditions which are required to be tested for the below result, mentioned in Lemma \ref{lem1} and \ref{lem2}. There has not been many known results exsistent for the $\zeta(1+ib)$ and hence we thought it would help in looking at such a computation.

\subsection{Background Results}
\label{lem1}\begin{lemma}
$f_n(x) = \frac{x^{-ib}}{x(x+n^a)}$ is $uniformly \hspace{2pt}convergent$ except the singularity at \{0\}
\end{lemma}
\begin{proofof} 1
Our proof shall be for $a=1$ cause we are looking at $\zeta(1+ib)$. The Uniform Norm Test in Section \ref{unm}  
cannot be applied because $x$ varies from 0 to $\infty$ in the integral and in the interval (0,1), $\frac{1}{x}$ turns out to be unbounded. Thus we proceed as follows:
$f_n(x)=\frac{e^{ib\ln(x)}}{x(x+n)} = \frac{\cos(b\ln(x))}{x(x+n)}-i\frac{\sin(b\ln(x))}{x(x+n)}$.
For a complex term $a_n(x)$+$ib_n(x)$ to be uniformly convergent both $a_n(x)$ and $b_n(x)$ should be uniformly convergent. We shall divide the above into 2 intervals \\
\begin{itemize}
\item  $x \in (1,\infty)$\\ $a_n(x) = \frac{\cos(b\ln(x))}{x(x+n)} \leq 1 $ for all $x,n \geq 0$\\
$\frac{\cos(b\ln(x))}{x(x+n)} \leq \frac{1}{x(x+n)}$ \\
$\lim_{n \rightarrow \infty} \frac{1}{x(x+n)} =0 \hspace{4pt} \forall x\geq 1$. Hence it is bounded. \\
Since $\lim_{n \rightarrow \infty} sup \Big\{ \frac{\cos(b\ln(x))}{x(x+n)}, \hspace{3pt} x\geq 1 \Big\}=\lim_{n \rightarrow \infty} \frac{1}{n+1}  =0$ \\
Thus $ \Big|\Big| \frac{\cos(b\ln(x))}{x(x+n)} \Big|\Big| \rightarrow 0$ and $ ||a_n - a|| \rightarrow 0$ \\
Similarly the same is true for $b_n(x)$ , hence $f_n(x)$ is uniformly convergent in the interval (1,$\infty$)
\item  $x \in (0,1)$ we rewrite the 
\begin{eqnarray}
a_n(x) =& \frac{\cos(b\ln(x))}{x(x+n)} \\
=&  \frac{\sum_{k=0}^{\infty} (1-x)^k \cos(b\ln(x))}{x+n}
\end{eqnarray} 
Thus as $\lim_{ \substack{n\rightarrow \infty ,\\ x\rightarrow 0}}$ the singularity has been removed and the $\int_0^1 \frac{cos(b\ln(x))}{x(x+n)}$ converges uniformly in the compact subset of any plane not containing 0. The above follows for $b_n(x)$ and thus $f_n(x)$ is uniformly convergent in the range (0,1]
\end{itemize}
Hence the uniform convergence of the $\int_0^{\infty} \frac{x^{-ib}}{x(x+n)}$ gives us the option of exchanging the summand and the integrand in Equation \ref{eqn2}
\end{proofof}\\
\label{lem2}\begin{lemma}
Both the functions (i)$f_{1n}(x)=\Big(\frac{-1}{2n+x+1}\Big)^k$,\hspace{2pt} (ii)$f_{2n}(x)=\Big(\frac{x^{-ib}}{(2n+x+1)(2n+x+2)}\Big)$ are uniformly convergent in the entire plane.
\end{lemma}
\begin{proofof} 2
(i) $f_{1n}(x)=(\frac{-1}{2n+x+1})^k \leq 1$ Hence $f_{1n}(x)$ is bounded. \\
$\lim_{n\rightarrow \infty}f_{1n}(x) =0 =f$  $ \forall\hspace{3pt} k \in [2,\infty)$ \hspace{3pt} $x \in \mcR^{+}$ + \{0\} and n \hspace{3pt} $\in \mcN + \{0\}$
\begin{eqnarray}
||f_{1n}(x) -f ||_{\mcR^{+}} &=& \Big|\Big| \Big(\frac{-1}{2n+x+1}\Big)^k \Big|\Big|_{\mcR^{+}} \\
&=& sup \Big\{ \Big| \Big(\frac{-1}{2n+x+1}\Big)^k \Big| : x\geq 0\Big\} \\
&=& (\frac{1}{2n})^k
\end{eqnarray}
We can infer that $\lim_{n \rightarrow \infty} (\frac{1}{2n})^{2k} \rightarrow 0$  $\forall k \in [2,\infty)$ \\
Hence $|| f_{1n}-f || \rightarrow 0$ . Thus $f_{1n(x)}$ is uniformly convergent to 0. \\
\begin{eqnarray*}
(ii)  f_{2n}(x) &=& \Big(\frac{x^{-ib}}{(2n+x+1)(2n+x+2)}\Big)  \\
&=& \Big(\frac{e^{-ib\ln(x)}}{(2n+x+1)(2n+x+2)}\Big)\\
&=& \Big(\frac{\cos(b\ln(x))}{(2n+x+1)(2n+x+2)}\Big)  - i \Big(\frac{\sin(b\ln(x))}{(2n+x+1)(2n+x+2)}\Big)
\end{eqnarray*}

$a_n(x)= \Big(\frac{\cos(b\ln(x))}{(2n+x+1)(2n+x+2)}\Big) \textless 1$. Hence $f_{2n}(x)$ is bounded. \\
$\lim_{n\rightarrow \infty} f_{2n}(x)=0 =f$  $\forall x \hspace{3 pt} \in \mcR^{+}$ + \{0\} and $n  \in \hspace{3 pt} \mcN + \{0\}$\\

\begin{eqnarray}
||a_n(x)-a||_{\mcR^{+}}=& \Big|\Big|\frac{\cos(b\ln(x))}{(2n+x+1)(2n+x+2)}\Big|\Big| \\
=& sup \Big\{ \Big| \frac{\cos(b\ln(x))}{(2n+x+1)(2n+x+2)} \Big|: x \geq 0 \Big\} \\
=& \frac{\cos(b\ln(x))}{(2n+1)(2n+2)} \leq \frac{1}{(2n+1)(2n+2)}
\end{eqnarray}
Since $\lim_{n\rightarrow \infty} \frac{1}{(2n+1)(2n+2)}=0$
Thus $|| a_n(x) - a|| \rightarrow 0$. The same follows for $b_n(x)$. 
Hence  $f_n(x)$ is uniformly convergent in the entire plane.
\end{proofof} \\
\subsection{Main Result}
The general expression for $\zeta(s)$ continued Re(s) $\textgreater$ 0 would be
\begin{eqnarray}
\zeta(s)=\frac{1}{1-2^{1-s}}\sum_{n=1}^{\infty} \frac{(-1)^n}{n^s}  s \neq 1
\end{eqnarray}
Let $s$=$a$+i$b$  $\forall$ $a$  $\textgreater$ 0 , $b$ $\neq$ 0

\label{equn1} \begin{equation} 
\zeta(a+ib) = \frac{1}{2^{1-a-ib}} \sum_{n=1}^{\infty} \frac{(-1)^n n^{-ib}}{n^a}
\end{equation}
Our aim or trick in handling the above expression is on the lines of \cite{TH07} where we would like to bring an integral within the summation. We would like to highlight that the following integral would help our expression 35.
\begin{eqnarray}
\int_{0}^{\infty} \frac{x^{-1} x^{-i(b/a)}}{x+n^a} dx &=& \frac{-\pi}{\sinh(\frac{ib\pi}{aO})} \frac{(n^{-a})^{-i(b/a)}}{n^a}  \\
\int_{0}^{\infty} \frac{x^{-1} x^{-ib}}{x+n} dx &=& \frac{-\pi}{\sinh(ib\pi)} \frac{(n^{-1})^{-ib}}{n} 
\end{eqnarray}
Using the above in Equation 35, we get an expression as follows 
\label{eqn2} \begin{eqnarray} 
\zeta(1+ib)= \frac{-1}{1-2^{-ib}} \frac{\sinh(ib\pi)}{\pi} \sum_{n=1}^{\infty} (-1)^n \int_{0}^{\infty} \frac{x^{-1} x^{-ib}}{x+n} dx
\end{eqnarray}

Using Lemma \ref{lem1}
\begin{eqnarray}
\zeta(1+ib) &=& \frac{-1}{1-2^{-ib}} \frac{\sinh(ib\pi)}{\pi} \sum_{n=1}^{\infty} (-1)^n \int_{0}^{\infty} \frac{x^{-1} x^{-ib}}{x+n} dx \\
&=& \frac{-1}{1-2^{-ib}} \frac{\sinh(ib\pi)}{\pi} \int_{0}^{\infty} \Big( \sum_{n=1}^{\infty} \frac{(-1)^n x^{-1} }{x+n}\Big) x^{-ib}  dx \label{eqn4}
\end{eqnarray}
We note that 
\begin{eqnarray}
\sum_{n=1}^{\infty} \frac{(-1)^n x^{-1} }{x+n} = \frac{1}{2} \Big( \Psi(\frac{x}{2}+1)-\Psi(\frac{x+1}{2})\Big )  
\end{eqnarray}
where the digamma function $\Psi(x)$= $\frac{\Gamma '(x)}{\Gamma (x)}$
\begin{eqnarray}
\zeta(1+ib)= \frac{-1}{1-2^{-ib}} \frac
{\sinh(ib\pi)}{2\pi} \int_{0}^{\infty} \Big( \Psi(\frac{x}{2}+1)-\Psi(\frac{x+1}{2})\Big ) x^{-ib} dx
\end{eqnarray}
We know the equation given in \cite{GR07} 
\begin{eqnarray}
\Psi(t+\alpha)-\Psi(\alpha)&=& \sum_{k=2}^{\infty} (-1)^k \zeta(k,\alpha)t^{k-1} \hspace{5 pt} \forall |t| \leq |\alpha|\\
\zeta(k,\alpha) &=& \sum_{n=0}^{\infty} \frac{1}{(n+\alpha)^k}  
\end{eqnarray}
With $\alpha=\frac{x+1}{2}$ and t=$\frac{1}{2}$. With x $\geq$ 0 . $\forall |t| \leq |\alpha| $
Hence
\begin{eqnarray}
\Psi(\frac{x}{2}+1)-\Psi(\frac{x+1}{2}) &=& \sum_{k=2}^{\infty} (-1)^k \zeta(k,\frac{x+1}{2}) (\frac{1}{2})^{k-1}\\
&=&\sum_{k=2}^{\infty} (-1)^k \sum_{n=0}^{\infty} \frac{2^k}{(2n+x+1)^k} 2^{1-k}\\
&=&\sum_{k=2}^{\infty} \sum_{n=0}^{\infty} \frac{2^k 2^{1-k} (-1)^k}{(2n+x+1)^k} \\
&=& 2\sum_{k=2}^{\infty} \sum_{n=0}^{\infty} (\frac{-1}{2n+x+1})^k 
\end{eqnarray}

Using the above in Equation \ref{eqn4} and Lemma \ref{lem2} we get the following 

\begin{eqnarray*}
\zeta(1+ib)&=& \frac{-1}{1-2^{-ib}} \frac{\sinh(ib\pi)}{2\pi} 2 \int_{0}^{\infty}\sum_{k=2}^{\infty} \sum_{n=0}^{\infty} \Big(\frac{-1}{2n+x+1}\Big)^k x^{-ib} dx \\
&=& \frac{-1}{1-2^{-ib}} \frac{\sinh(ib\pi)}{\pi} \lim_{k \rightarrow \infty} \int_{0}^{\infty} \sum_{n=0}^{k} \frac{x^{-ib}}{(2n+x+1)(2n+x+2)} dx
\end{eqnarray*}
With Lemma \ref{lem2} we can simplify the above to
\begin{multline}
\zeta(1+ib)= \frac{-1}{1-2^{-ib}} \frac{\sinh(ib\pi)}{\pi} \sum_{n=0}^{\infty} (2n+1)^{-ib} \Big( -i\pi csch(i\pi b) \\+i \pi csch(i\pi b) ((2n+1)^{ib} (2n+2)^{-ib}) \Big)\\
\end{multline}
\begin{eqnarray}
\zeta(1+ib)&=& \frac{-1}{1-2^{-ib}} \sum_{n=0}^{\infty} ( -(2n+1)^{-ib} + (2n+2)^{-ib})\\
&=&  \frac{-1}{1-2^{-ib}} \sum_{n=1}^{\infty} \frac{(-1)^n}{n^{ib}}
\end{eqnarray}
Basically the above result stems from the first step but theres a subtle factor of $n$ vanishes in the denominator which emphasis that the above is the actual formula for $\zeta(1+ib)$. The rigorous need for proofs of uniform convergence show the innate siginificance of convergence when related to $zeta(s)$. We could validate few known results like $\lim_{n\rightarrow \infty} (s-1)\zeta(s)=1$ , $\lim_{b\rightarrow 0} \zeta(1+ib)\rightarrow \infty$ , there is also the mention that the other zeros would lie on $b=\frac{2k\pi}{ \log2}$ for any $k \neq 0\in \mcI$ ( which can be got by simply relating the analytic continuation of Dirichelet Series and the Alternating Zeta function). Hence we could find more identities with regard to the zeta function on the Line $s=1$ as shown in \cite{S03}

\section{Conclusion}
This article gives a look into the evaluation of the Riemann Zeta function with odd arguments and at the strip $\Re(s)=1$ such that we could achieve more results and analyze the regions with better knowledge. The advantages of each formulation have been mentioned in each section. We might compromised a little on the accuracy but our main goal has been to reduce complexity(and providing exponential convergence) so that the outcome would be almost accurate but evaluated faster with lesser parameters and fewer huge calculations. We also have given a result for the $\zeta(1+ib)$ so that better analysis could be made with the formula.

\newpage

\newcommand{\etalchar}[1]{$^{#1}$}

\end{document}